\documentclass[11pt]{amsart}

\theoremstyle{plain}
\newtheorem{theorem}                {Theorem}      [section]
\newtheorem{proposition}  [theorem]  {Proposition}
\newtheorem{corollary}    [theorem]  {Corollary}
\newtheorem{lemma}        [theorem]  {Lemma}

\theoremstyle{definition}

\newtheorem{remark}       [theorem]  {Remark}

 \DeclareMathOperator{\id}{I}

\DeclareMathOperator{\Span}{span}
 
\DeclareMathOperator{\cst}{constant}
 
 \DeclareMathOperator{\im}{Im}

\numberwithin{equation}{section}

\begin{document}

\title[Surfaces with parallel mean curvature vector]
{Surfaces with parallel mean curvature vector in complex space
forms}

\author{Dorel Fetcu}

\address{Department of Mathematics\\
"Gh. Asachi" Technical University of Iasi\\
Bd. Carol I no. 11 \\
700506 Iasi, Romania} \email{dfetcu@math.tuiasi.ro}

\curraddr{IMPA\\ Estrada Dona Castorina 110, 22460-320 Rio de
Janeiro, Brazil} \email{dorel@impa.br}

\thanks{The author was supported by a Post-Doctoral Fellowship "P\'os-Doutorado J\'unior
(PDJ)" offered by CNPq, Brazil.}

\begin{abstract} We consider a quadratic form defined on the
surfaces with parallel mean curvature vector of an any dimensional
complex space form and prove that its $(2,0)$-part is holomorphic.
When the complex dimension of the ambient space is equal to $2$ we
define a second quadratic form with the same property and then
determine those surfaces with parallel mean curvature vector on
which the $(2,0)$-parts of both of them vanish. We also provide a
reduction of codimension theorem and prove a non-existence result
for $2$-spheres with parallel mean curvature vector.
\end{abstract}

\date{}



\maketitle

\section{Introduction}

Almost sixty years ago H. Hopf was the first to use a quadratic
form in order to study surfaces immersed in a $3$-dimensional
Euclidean space. He proved, in 1951, that any such surface which
is homeomorphic to a sphere and has constant mean curvature is
actually isometric to a round sphere (see \cite{HH}). This result
was extended by S.-S. Chern to surfaces immersed in
$3$-dimensional space forms (see \cite{SSC}) and by U. Abresch and
H. Rosenberg to surfaces in simply connected, homogeneous
$3$-dimensional Riemannian manifolds, whose group of isometries
has dimension $4$ (see \cite{AR, AR2}). Very recently, H. Alencar,
M. do Carmo and R. Tribuzy have made the next step by obtaining
Hopf-type results in spaces with dimension higher than $3$, namely
in product spaces $M^n(c)\times\mathbb{R}$, where $M^n(c)$ is a
simply connected $n$-dimensional space form with constant sectional curvature
$c\neq 0$ (see \cite{AdCT2, AdCT}). They have considered the case
of surfaces with parallel mean curvature vector, as a natural
generalization of those with constant mean curvature in a
$3$-dimensional ambient space. We also have to mention a recent
paper of F. Torralbo and F. Urbano, which is devoted to the study
of surfaces with parallel mean curvature vector in
$\mathbb{S}^2\times\mathbb{S}^2$ and
$\mathbb{H}^2\times\mathbb{H}^2$.

Minimal surfaces and surfaces with parallel mean curvature vector
in complex space forms have been also a well studied subject in
the last two decades (see, for example,
\cite{BJRW,CO,CW,EGT,H,KZ,O,O2,YO}). In all these papers the
K\"ahler angle proved to play a decisive role in understanding of
the geometry of immersed surfaces in a complex space form, and, in
several of them, important results were obtained when this angle
was supposed to be constant (see \cite{BJRW,O,YO}).

The main goal of our paper is to obtain some characterization
results concerning surfaces with parallel mean curvature vector in
complex space forms by using as a principal tool holomorphic
quadratic forms defined on these surfaces. The paper is organized
as follows. In Section $2$ we introduce a quadratic form $Q$ on
surfaces of an arbitrary complex space form and prove that its
$(2,0)$-part is holomorphic when the mean curvature vector of the
surface is parallel. In Section $3$ we work in the complex space
forms with complex dimension equal to $2$ and find another
quadratic form $Q'$ with holomorphic $(2,0)$-part. Then we
determine surfaces with parallel mean curvature vector on which
both $(2,0)$-part of $Q$ and $(2,0)$-part of $Q'$ vanish. As a
by-product we reobtain a result in \cite{H}. More precisely, we
prove that a $2$-sphere can be immersed as a surface with parallel
mean curvature vector only in a flat complex space form and it is
a round sphere in a hyperplane in $\mathbb{C}^2$. In Section $4$
we deal with surfaces in $\mathbb{C}^n$ with parallel mean
curvature vector, and we prove that the $(2,0)$-part of $Q$
vanishes on such a surface if and only if it is pseudo-umbilical.
The main result of Section $5$ is a reduction theorem, which
states that a surface in a complex space form, with parallel mean
curvature vector, either is totally real and pseudo-umbilical or
it is not pseudo-umbilical and lies in a complex space form with
complex dimension less or equal to $5$. The last Section is
devoted to the study of the $2$-spheres with parallel
mean curvature vector and constant K\"ahler angle. We prove that
there are no non-pseudo-umbilical such spheres in a complex space
form with constant holomorphic sectional curvature $\rho\neq 0$.

\noindent \textbf{Acknowledgements.} The author wants to thank
Professor Harold Rosenberg for suggesting this subject, useful
comments and discussions and constant encouragement.

\section{A quadratic form}

Let $\Sigma^2$ be an immersed surface in $N^{n}(\rho)$, where $N$
is a complex space form with complex dimension $n$, complex
structure $(J,\langle,\rangle)$, and with constant holomorphic
sectional curvature $\rho$; which is $\mathbb{C}P^n(\rho)$,
$\mathbb{C}^n$ or $\mathbb{C}H^n(\rho)$, as $\rho>0$, $\rho=0$ and
$\rho<0$, respectively. Let us define a quadratic form $Q$ on
$\Sigma^2$ by
$$
Q(X,Y)=8|H|^2\langle\sigma(X,Y),H\rangle+3\rho\langle JX,
H\rangle\langle JY, H\rangle,
$$
where $\sigma$ is the second fundamental form of the surface and
$H$ is its mean curvature vector field. Assume that $H$ is
parallel in the normal bundle of $\Sigma^2$, i.e.
$\nabla^{\perp}H=0$, the normal connection $\nabla^{\perp}$ being
defined by the equation of Weingarten
$$
\nabla^{N}_XV=-A_VX+\nabla^{\perp}_XV,
$$
for any vector field $X$ tangent to $\Sigma^2$ and any vector
field $V$ normal to the surface, where $\nabla^{N}$ is the
Levi-Civita connection on $N$ and $A$ is the shape operator.

We shall prove that the $(2,0)$-part of $Q$ is holomorphic. In
order to do that, let us first consider the isothermal coordinates
$(u,v)$ on $\Sigma^2$. Then $ds^2=\lambda^2(du^2+dv^2)$ and define
$z=u+iv$, $\bar z=u-iv$, $dz=\frac{1}{\sqrt{2}}(du+idv)$, $d\bar
z=\frac{1}{\sqrt{2}}(du-idv)$ and
$$
Z=\frac{1}{\sqrt{2}}\Big(\frac{\partial}{\partial
u}-i\frac{\partial}{\partial v}\Big),\quad \bar
Z=\frac{1}{\sqrt{2}}\Big(\frac{\partial}{\partial
u}+i\frac{\partial}{\partial v}\Big).
$$
We also have $\langle Z,\bar
Z\rangle=\langle\frac{\partial}{\partial
u},\frac{\partial}{\partial
u}\rangle=\langle\frac{\partial}{\partial
v},\frac{\partial}{\partial v}\rangle=\lambda^2$.

In the following we shall calculate
$$
\bar Z(Q(Z,Z))=\bar
Z(8|H|^2\langle\sigma(Z,Z),H\rangle+3\rho\langle JZ, H\rangle^2).
$$

First, we get
$$
\begin{array}{ll}
\bar Z(\langle\sigma(Z,Z),H\rangle)&=\langle\nabla^N_{\bar
Z}\sigma(Z,Z),H\rangle+\langle\sigma(Z,Z),\nabla^N_{\bar Z}H\rangle\\
\\&=\langle\nabla^{\perp}_{\bar Z}\sigma(Z,Z),H\rangle+\langle\sigma(Z,Z),\nabla^{\perp}_{\bar
Z}H\rangle\\ \\&=\langle(\nabla^{\perp}_{\bar
Z}\sigma)(Z,Z),H\rangle+\langle\sigma(Z,Z),\nabla^{\perp}_{\bar
Z}H\rangle,
\end{array}
$$
where we have used that
$$
(\nabla^{\perp}_{\bar Z}\sigma)(Z,Z)=\nabla^{\perp}_{\bar
Z}\sigma(Z,Z)-2\sigma(\nabla_{\bar Z}Z,Z)=\nabla^{\perp}_{\bar
Z}\sigma(Z,Z)
$$
since, from the definition of the connection $\nabla$ on the
surface, we easily get $\nabla_{\bar Z}Z=0$.

Now, from the Codazzi equation, we obtain
\begin{equation}\label{eq:1}
\begin{array}{lll}
\bar Z(\langle
\sigma(Z,Z),H\rangle)&=&\langle(\nabla^{\perp}_{Z}\sigma)(\bar
Z,Z),H\rangle+\langle (R^N(\bar
Z,Z)Z)^{\perp},H\rangle\\ \\&&+\langle\sigma(Z,Z),\nabla^{\perp}_{\bar Z}H\rangle\\
\\&=&\langle(\nabla^{\perp}_{Z}\sigma)(\bar Z,Z),H\rangle+\langle R^N(\bar
Z,Z)Z,H\rangle+\langle\sigma(Z,Z),\nabla^{\perp}_{\bar Z}H\rangle.
\end{array}
\end{equation}

From the expression of the curvature tensor field of $N$
$$
\begin{array}{lcl}
R^N(U,V)W&=&\frac{\rho}{4}\{\langle V,W\rangle U-\langle U,W\rangle V+\langle JV,W\rangle JU-\langle JU,W\rangle JV\\
\\&&+2\langle JV,U\rangle JW\},
\end{array}
$$
it follows
\begin{equation}\label{eq:2}
\langle R^N(\bar Z,Z)Z,H\rangle=\frac{3\rho}{4}\langle\bar
Z,JZ\rangle\langle H,JZ\rangle.
\end{equation}

We also have the following

\begin{lemma}
\begin{equation}\label{eq:3}
\langle(\nabla^{\perp}_{Z}\sigma)(\bar Z,Z),H\rangle=\langle\bar
Z,Z\rangle\langle\nabla^{\perp}_ZH,H\rangle.
\end{equation}
\end{lemma}

\begin{proof} By using the definition of $(\nabla^{\perp}_{Z}\sigma)(\bar
Z,Z)$ one obtains
$$
(\nabla^{\perp}_{Z}\sigma)(\bar Z,Z)=\nabla^{\perp}_{Z}\sigma(\bar
Z,Z)-\sigma(\nabla_Z\bar Z,Z)-\sigma(\bar
Z,\nabla_ZZ)=\nabla^{\perp}_{Z}\sigma(\bar Z,Z)-\sigma(\bar
Z,\nabla_ZZ)
$$
since $\nabla_Z\bar Z=0$.

Next, let us consider the unit vector fields $e_1$ and $e_2$
corresponding to $\frac{\partial}{\partial u}$ and
$\frac{\partial}{\partial v}$, respectively, and
$E=\frac{1}{\sqrt{2}}(e_1-ie_2)$. Then we have $Z=\lambda E$ and
$$
\sigma(\bar
Z,Z)=\frac{\lambda^2}{2}\sigma(e_1-ie_2,e_1+ie_2)=\frac{\lambda^2}{2}(\sigma(e_1,e_1)+\sigma(e_2,e_2))=\langle\bar Z,Z\rangle H.
$$

Since $\nabla_ZZ$ is tangent it follows that $\nabla_ZZ=a Z+b\bar
Z$ and then $0=\langle\nabla_ZZ,Z\rangle=b\lambda^2$, where we
have used the fact that $\langle Z,Z\rangle=0$, and
$a=\frac{1}{\lambda^2}\langle\nabla_ZZ,\bar Z\rangle$.

In conclusion
$$
\begin{array}{lll}
\langle(\nabla^{\perp}_{Z}\sigma)(\bar
Z,Z),H\rangle&=&\langle\nabla^N_Z(\langle\bar Z,Z\rangle H),H\rangle-\langle\nabla_ZZ,\bar Z\rangle\langle H,H\rangle\\
\\&=&\langle\nabla_Z\bar Z,Z\rangle\langle H,H\rangle+\langle\nabla_ZZ,\bar Z\rangle\langle H,H\rangle\\
\\&&+\langle\bar Z,Z\rangle\langle\nabla^{\perp}_ZH,H\rangle-\langle\nabla_ZZ,\bar Z\rangle\langle H,H\rangle\\
\\&=&\langle\bar Z,Z\rangle\langle\nabla^{\perp}_ZH,H\rangle.
\end{array}
$$
\end{proof}

\begin{lemma}
\begin{equation}\label{eq:4}
\bar Z(\langle JZ,H\rangle^2)=2\langle
JZ,H\rangle\langle(JZ)^{\perp},\nabla^{\perp}_{\bar
Z}H\rangle-2|H|^2\langle\bar Z,JZ\rangle\langle JZ,H\rangle
\end{equation}
\end{lemma}

\begin{proof}
From the definitions of the K\"ahler structure and of the
Levi-Civita connection we have
$$
\begin{array}{lll}
\bar Z(\langle JZ,H\rangle^2)&=&2\langle
JZ,H\rangle\{\langle\nabla^N_{\bar
Z}JZ,H\rangle+\langle JZ,\nabla^N_{\bar Z}H\rangle\}\\
\\&=&2\langle JZ,H\rangle\{\langle\bar
Z,Z\rangle\langle JH,H\rangle-\langle(JZ)^{\top},A_H\bar Z\rangle\\
\\&&+\langle(JZ)^{\perp},\nabla^{\perp}_{\bar Z}H\rangle\}\\ \\
&=&2\langle JZ,H\rangle\{\langle(JZ)^{\perp},\nabla^{\perp}_{\bar
Z}H\rangle-\langle\sigma((JZ)^{\top},\bar Z),H\rangle\}\\
\\&=&2\langle JZ,H\rangle\{\langle(JZ)^{\perp},\nabla^{\perp}_{\bar Z}H\rangle-\langle JZ,\bar
Z\rangle|H|^2\},
\end{array}
$$
where we have used $\nabla^N_{\bar Z}Z=\sigma(\bar
Z,Z)=\langle\bar Z,Z\rangle H$, as we have seen in the proof of
the previous Lemma, and $(JZ)^{\top}=\frac{1}{\lambda^2}\langle
JZ,\bar Z\rangle Z$, that can be easily checked.
\end{proof}

By replacing \eqref{eq:2}, \eqref{eq:3} and \eqref{eq:4} into
\eqref{eq:1} we obtain that $\bar Z(Q(Z,Z))$ vanishes and then we
come to the conclusion that

\begin{proposition}
If $\Sigma^2$ is an immersed surface in a complex space form
$N^{n}(\rho)$, with parallel mean curvature vector field, then the
$(2,0)$-part of the quadratic form $Q$, defined on $\Sigma^2$ by
$$
Q(X,Y)=8|H|^2\langle\sigma(X,Y),H\rangle+3\rho\langle JX,
H\rangle\langle JY, H\rangle,
$$
is holomorphic.
\end{proposition}

\section{Quadratic forms and $2$-Spheres in $2$-dimensional complex space forms}

In this section we shall define a new quadratic form on a surface
$\Sigma^2$ immersed in a complex space form $N^2(\rho)$, with
parallel mean curvature vector field $H\neq 0$, and prove that its
$(2,0)$-part is holomorphic. Then, by using these two quadratic
forms, we shall classify the $2$-spheres with nonzero parallel
mean curvature vector.

\subsection{Another quadratic form}

Let us consider an oriented orthonormal local frame $\{\widetilde
e_1,\widetilde e_2\}$ on the surface and denote by $\theta$ the
K\"ahler angle function defined by
$$
\langle J\widetilde e_1,\widetilde e_2\rangle=\cos\theta.
$$
The immersion $x:\Sigma^2\rightarrow N$ is said to be holomorphic
if $\cos\theta=1$, anti-holomorphic if $\cos\theta=-1$, and
totally real if $\cos\theta=0$. In the following we shall assume
that $x$ is neither holomorphic or anti-holomorphic.

Next, we take $e_3=-\frac{H}{|H|}$ and let $e_4$ be the unique
unit normal vector field orthogonal to $e_3$ compatible with the
orientation of $\Sigma^2$ in $N$. Since $e_3$ is parallel in the
normal bundle so is $e_4$, and, as the K\"ahler angle is
independent of the choice of the orthonormal frame on the surface
(see, for example, \cite{CW}), we have
\begin{equation}\label{eq0'}
\langle Je_4,e_3\rangle=\cos\theta.
\end{equation}

Now, we can consider the vector fields
$$
e_1=\cot\theta e_3-\frac{1}{\sin\theta}Je_4,\quad
e_2=\frac{1}{\sin\theta}Je_3+\cot\theta e_4
$$
tangent to the surface and obtain an orthonormal frame field
$\{e_1,e_2,e_3,e_4\}$ adapted to $\Sigma^2$ in $N$.

We define a quadratic form $Q'$ on $\Sigma^2$ by
$$
Q'(X,Y)=8i|H|\langle\sigma(X,Y),e_4\rangle+3\rho\langle
JX,e_4\rangle\langle JY,e_4\rangle
$$
and again consider the isothermal coordinates $(u,v)$ on
$\Sigma^2$ and the tangent complex vector fields $Z$ and $\bar Z$.
In the same way as in the case of $Q$, using the Codazzi equation,
the fact that $H$ and $e_4$ are parallel and the expression of the
curvature vector field of $N$, we get
\begin{equation}\label{eq1'}
\bar Z(\langle\sigma(Z,Z),e_4\rangle)=\frac{3\rho}{4}\langle\bar
Z,JZ\rangle\langle JZ,e_4\rangle.
\end{equation}
On the other hand, we have
$$
\begin{array}{lll}
\bar Z(\langle JZ,e_4\rangle^2)&=&2\langle
JZ,e_4\rangle\{\langle\nabla^N_{\bar
Z}JZ,e_4\rangle+\langle JZ,\nabla^N_{\bar Z}e_4\rangle\}\\
\\&=&2\langle JZ,e_4\rangle\{\langle\bar
Z,Z\rangle\langle JH,e_4\rangle-\langle(JZ)^{\top},A_{e_4}\bar Z\rangle\}\\
\\&=&-2|H|\langle JZ,e_4\rangle\langle\bar Z,Z\rangle\langle Je_3,e_4\rangle-2\langle JZ,e_4\rangle\langle\sigma((JZ)^{\top},\bar Z),e_4\rangle\\
\\&=&2|H|\langle JZ,e_4\rangle\langle\bar Z,Z\rangle\cos\theta-2\langle JZ,e_4\rangle\langle JZ,\bar Z\rangle\langle H,e_4\rangle\\
\\&=&2|H|\langle JZ,e_4\rangle\langle \bar Z,Z\rangle\cos\theta,
\end{array}
$$
where we have used $\nabla^N_{\bar Z}Z=\sigma(\bar
Z,Z)=\langle\bar Z,Z\rangle H$,
$(JZ)^{\top}=\frac{1}{\lambda^2}\langle JZ,\bar Z\rangle Z$ and
\eqref{eq0'}. But $\langle\bar Z,JZ\rangle=-i\langle\bar
Z,Z\rangle\langle e_1,Je_2\rangle=i\langle\bar
Z,Z\rangle\cos\theta$, and therefore
\begin{equation}\label{eq2'}
\bar Z(\langle JZ,e_4\rangle^2)=-2i|H|\langle\bar
Z,JZ\rangle\langle JZ,e_4\rangle.
\end{equation}
Hence, from \eqref{eq1'} and \eqref{eq2'}, one obtains $\bar
Z(Q'(Z,Z))=0$, which means that the $(2,0)$-part of the quadratic
form $Q'$ is holomorphic.

\subsection{$2$-Spheres in $2$-dimensional complex space forms}

In order to classify the $2$-spheres in $2$-dimensional complex
space forms, we shall need a result of T. Ogata in \cite{O}, which
we will briefly recall in the following (see also \cite{H} and
\cite{KZ}). Consider a surface $\Sigma^2$ isometrically immersed
in a complex space form $N^2(\rho)$, with parallel mean curvature
vector field $H\neq 0$. Using the frame field on $N^2(\rho)$
adapted to $\Sigma^2$, defined above, and considering isothermal
coordinates $(u,v)$ on the surface, Ogata proved that there exist
complex-valued functions $a$ and $c$ on $\Sigma^2$ such that
$\theta$, $\lambda$, $a$ and $c$ satisfy
\begin{equation}\label{eq1s}
\begin{cases}
\frac{\partial\theta}{\partial
z}=\lambda(a+b)\\\frac{\partial\lambda}{\partial\bar
z}=-|\lambda|^2(\bar a-b)\cot\theta\\\frac{\partial
a}{\partial\bar
z}=\bar\lambda\Big(2|a|^2-2ab+\frac{3\rho\sin^2\theta}{8}\Big)\cot\theta\\\frac{\partial
c}{\partial
z}=2\lambda(a-b)c\cot\theta\\|c|^2=|a|^2+\frac{\rho(3\sin^2\theta-2)}{8}
\end{cases}
\end{equation}
where $z=u+iv$ and $|H|=2b$; and also the converse: if $\rho$ is a
real constant, $b$ a positive constant, $\Sigma^2$ a
$2$-dimensional Riemannian manifold, and there exist some
functions $\theta$, $a$ and $c$ on $\Sigma^2$ satisfying
\eqref{eq1s}, then there is an isometric immersion of $\Sigma^2$
into $N^2(\rho)$ with parallel mean curvature vector field of
length equal to $2b$ and with the K\"ahler angle $\theta$. The
second fundamental form of $\Sigma^2$ in $N$ w.r.t.
$\{e_1,e_2,e_3,e_4\}$ is given by
\begin{footnotesize}
$$
\sigma^3=\left(\begin{array}{cc}-2b-\Re(\bar a+c)&-\Im(\bar a+c)\\
\\-\Im(\bar a+c)&-2b+\Re(\bar
a+c)\end{array}\right)\quad\textnormal{and}\quad
\sigma^4=\left(\begin{array}{cc}\Im(\bar a-c)&-\Re(\bar a-c)\\
\\-\Re(\bar a-c)&-\Im(\bar a-c)\end{array}\right)
$$
\end{footnotesize}
and the Gaussian curvature of $\Sigma^2$ is
$K=4b^2-4|c|^2+\frac{\rho}{2}$ (see also \cite{H}).

Assume now that the $(2,0)$-part of $Q$ and the $(2,0)$-part of
$Q'$ vanish on the surface $\Sigma^2$. It follows, from the
expression of the second fundamental form, that $\bar
c+a\in\mathbb{R}$, $\bar c-a\in\mathbb{R}$ and
$$
32b(\bar c+a)-3\rho\sin^2\theta=0,\quad 32b(\bar
c-a)+3\rho\sin^2\theta=0.
$$
Therefore $c=0$ and $a=\frac{3\rho\sin^2\theta}{32b}$ and, from
the fifth equation of \eqref{eq1s}, it follows
\begin{equation}\label{eq2s}
9\rho^2\sin^4\theta+128\rho b^2(3\sin^2\theta-2)=0.
\end{equation}
We have to split the study of this equation in two cases. First,
if $\rho=0$ then the above equation holds and $a=0$. Next, if
$\rho\neq 0$, we get that function $\theta$ is a constant. This,
together with the first equation of \eqref{eq1s}, lead to
$a=\frac{3\rho\sin^2\theta}{32b}=-b$. By replacing in equation
\eqref{eq2s} we obtain $\rho=-12b^2$ and then
$\sin^2\theta=\frac{8}{9}$. We note that in both cases the
Gaussian curvature of $\Sigma^2$ is given by
$K=4b^2+\frac{\rho}{2}=\cst$ (see \cite{H}). Thus, by using
Theorem 1.1 in \cite{H}, we have just proved that

\begin{theorem} If the $(2,0)$-part of $Q$ and the $(2,0)$-part of $Q'$ vanish on a
surface $\Sigma^2$ isometrically immersed in a complex space form
$N^2(\rho)$, with parallel mean curvature vector field of length
$2b>0$, then either
\begin{enumerate}
\item $N^2(\rho)=\mathbb{C}H^2(-12b^2)$ and $\Sigma^2$ is the
slant surface in \cite{C} (Theorem 3(2));

\item $N^2(\rho)=\mathbb{C}^2$ and $\Sigma^2$ is a part of a round
sphere in a hyperplane in $\mathbb{C}^2$.
\end{enumerate}
\end{theorem}

Since the Gaussian curvature $K$ is nonnegative only in the second
case of the Theorem, we have also reobtained the following result
of S. Hirakawa in \cite{H}.

\begin{corollary}\label{th_sphere_2} If $\mathbb{S}^2$ is an isometrically immersed sphere in a
$2$-dimensional complex space form, with nonzero parallel mean
curvature vector, then it is a round sphere in a hyperplane in
$\mathbb{C}^2$.
\end{corollary}

\section{A remark on the $2$-spheres in $\mathbb{C}^n$}

\begin{proposition}\label{p_spheres_n} Let $\Sigma^2$ be an isometrically immersed
surface in $\mathbb{C}^n$, with nonzero parallel mean curvature
vector. Then the $(2,0)$-part of the quadratic form $Q$ vanishes
on $\Sigma^2$ if and only if the surface is pseudo-umbilical, i.e.
$A_H=|H|^2\id$.
\end{proposition}

\begin{proof} It can be easily seen that if $\Sigma^2$ is pseudo-umbilical
then the $(2,0)$-part of $Q$ vanishes and, therefore, we have to
prove only the necessity.

From $Q(Z,Z)=\frac{\langle Z,\bar
Z\rangle^2}{2}Q(e_1-ie_2,e_1-ie_2)=0$ it follows
$$
\langle\sigma(e_1,e_1)-\sigma(e_2,e_2),H\rangle=0
$$
and
$$
\langle\sigma(e_1,e_2),H\rangle=0.
$$
But, since
$\langle\sigma(e_1,e_1)+\sigma(e_2,e_2),H\rangle=2|H|^2$, we
obtain, for each $i\in\{1,2\}$,
$$
\langle A_He_i,e_i\rangle=\langle\sigma(e_i,e_i),H\rangle=|H|^2.
$$
Therefore $A_H=|H|^2\id$, i.e. $\Sigma^2$ is pseudo-umbilical.
\end{proof}

S.-T. Yau proved (Theorem 4 in \cite{Y}) that if $\Sigma^2$ is a
surface with parallel mean curvature vector $H$ in a manifold $N$
with constant sectional curvature, then either $\Sigma^2$ is a
minimal surface of an umbilical hypersurface of $N$ or $\Sigma^2$
lies in a $3$-dimensional umbilical submanifold of $N$ with
constant mean curvature, as $H$ is an umbilical direction or the
second fundamental form of $\Sigma^2$ can be diagonalized
simultaneously. We note that, in the first case, the mean
curvature vector field of $\Sigma^2$ in $\mathbb{C}^n$ is
orthogonal to the hypersurface.

Applying this result, together with Proposition \ref{p_spheres_n},
to the $2$-spheres in $\mathbb{C}^n$, and using the Gauss equation
of a hypersurface in $\mathbb{C}^n$, we get

\begin{proposition}\label{th_sphere_n} If $\mathbb{S}^2$ is an isometrically
immersed sphere in $\mathbb{C}^n$, with nonzero parallel mean
curvature vector field $H$, then it is a minimal surface of a
hypersphere $\mathbb{S}^{2n-1}(|H|)\subset\mathbb{C}^n$.
\end{proposition}

\section{Reduction of the codimension}

Let $x:\Sigma^2\rightarrow N^n(\rho)$, $n\geq 3$, $\rho\neq 0$, be
an isometric immersion of a surface $\Sigma^2$ in a complex space
form, with parallel mean curvature vector field $H\neq 0$.

\begin{lemma}\label{lemma_com} For any vector $V$ normal to $\Sigma^2$, which is also orthogonal to
$JT\Sigma^2$ and to $JH$, we have $[A_H,A_V]=0$, i.e. $A_H$
commutes with $A_V$.
\end{lemma}

\begin{proof}  The statement follows easily, from the Ricci equation
$$
\langle
R^{\perp}(X,Y)H,V\rangle=\langle[A_H,A_V]X,Y\rangle+\langle
R^N(X,Y)H,V\rangle,
$$
since
$$
\begin{array}{lll}
\langle R^N(X,Y)H,V\rangle&=&\frac{\rho}{4}\{\langle JY,H\rangle\langle JX,V\rangle-\langle JX,H\rangle\langle JY,V\rangle\\
\\&&+2\langle JY,X\rangle\langle JH,V\rangle\}\\ \\&=&0
\end{array}
$$
and $R^{\perp}(X,Y)H=0$.
\end{proof}

\begin{remark} If $n=3$ and $H\perp JT\Sigma^2$ do not hold simultaneously,
then there exists at least one normal vector $V$ as in Lemma
\ref{lemma_com}. This can be proved by using the basis of the
tangent space $TN$ along $\Sigma^2$ defined in \cite{O2}, which
construction we shall briefly explain in the following. Let us
consider a local orthonormal frame $\{e_1,e_2\}$ of vector fields
tangent to $\Sigma^2$. Since we have assumed that $H\neq 0$ it
follows that $\Sigma^2$ is not holomorphic or antiholomorphic,
which means that $\cos^2\theta=1$ only at isolated points, and we
shall work in the open dense set of points where $\cos^2\theta\neq
1$, where $\theta$ is the K\"ahler angle function. The next step
is to define two normal vectors by
$$
e_3=-\cot\theta
e_1-\frac{1}{\sin\theta}Je_2\quad\textnormal{and}\quad
e_4=\frac{1}{\sin\theta}Je_1-\cot\theta e_2
$$
and now we have an orthonormal basis $\{e_1,e_2,e_3,e_4\}$ of
$\Span\{e_1,e_2,Je_1,Je_2\}$. Moreover, we can set
$$
\widetilde
e_1=\cos\Big(\frac{\theta}{2}\Big)e_1+\sin\Big(\frac{\theta}{2}\Big)e_3,\quad\widetilde
e_2=\cos\Big(\frac{\theta}{2}\Big)e_2+\sin\Big(\frac{\theta}{2}\Big)e_4
$$
$$
\widetilde
e_3=\sin\Big(\frac{\theta}{2}\Big)e_1-\cos\Big(\frac{\theta}{2}\Big)e_3,\quad\widetilde
e_4=-\sin\Big(\frac{\theta}{2}\Big)e_2+\cos\Big(\frac{\theta}{2}\Big)e_4
$$
and obtain a $J$-canonical basis of $\Span\{e_1,e_2,Je_1,Je_2\}$,
i.e. $J\widetilde e_{2i-1}=\widetilde e_{2i}$. Finally, let us
consider a $J$-basis of $TN$ along $\Sigma^2$, of the form
$\{\widetilde e_1, \widetilde e_2, \widetilde e_3, \widetilde e_4,
\widetilde e_5, \widetilde e_6=J\widetilde e_5,\ldots,\widetilde
e_{2n-1},\widetilde e_{2n}=J\widetilde e_{2n-1}\}$. Now, three
situations can occur:
\begin{enumerate}

\item $H\in(JT\Sigma^2)^{\perp}$, and then $\widetilde e_5\perp
JT\Sigma^2$ and $\widetilde e_5\perp JH$, where we have denoted by
$(JT\Sigma^2)^{\perp}=\{(JX)^{\perp}:X\ \textnormal{tangent to}\
\Sigma^2\}$;

\item $H\perp JT\Sigma^2$, and then, if we choose $\widetilde
e_5=H$ and $\widetilde e_6=JH$, we have $\widetilde e_7\perp
JT\Sigma^2$ and $\widetilde e_7\perp JH$ (obviously, this case can
occur only if $n>3$);

\item $H\notin(JT\Sigma^2)^{\perp}$ and $H$ is not orthogonal to
$JT\Sigma^2$. In this case we may consider the vector $u$, the
projection of $H$ on the complementary space of
$(JT\Sigma^2)^{\perp}$ in $TN$ (along $\Sigma^2$) and set
$\widetilde e_5=\frac{u}{|u|}$. It follows that $\widetilde
e_5\perp JT\Sigma^2$ and $\widetilde e_5\perp JH$.
\end{enumerate}

If $n=3$ and $H\perp JT\Sigma^2$ it is easy to see that
$$
\langle R^N(X,Y)H,e_3\rangle=\langle R^N(X,Y)H,e_4\rangle=0
$$
for any vector fields $X$ and $Y$ tangent to $\Sigma^2$, and then
that $A_H$ commutes with $A_{e_3}$ and $A_{e_4}$.
\end{remark}

Conclusively, we get the following

\begin{corollary}\label{lemma_split} Either $H$ is an
umbilical direction or there exists a basis that diagonalizes
simultaneously $A_H$ and $A_V$, for all normal vectors satisfying
$V\perp JH$, if $n=3$ and $H\perp JT\Sigma^2$, or the conditions
in Lemma \ref{lemma_com}, otherwise.
\end{corollary}

\begin{lemma}\label{lemma_parallel} Assume that $H$ is nowhere an umbilical direction.
Then there exists a parallel subbundle of the normal bundle which
contains the image of the second fundamental form $\sigma$ and has
dimension less or equal to $8$.
\end{lemma}

\begin{proof} We consider the following subbundle $L$ of the normal bundle
$$
L=\Span\{\im\sigma\cup(J\im\sigma)^{\perp}\cup(JT\Sigma^2)^{\perp}\},
$$
and we will show that $L$ is parallel.

First, we shall prove that, if $V$ is orthogonal to $L$, then
$\nabla^{\perp}_{e_i}V$ is orthogonal to $JT\Sigma^2$ and to $JH$,
where $\{e_1,e_2\}$ is a frame w.r.t. which we have
$\langle\sigma(e_1,e_2),V\rangle=\langle\sigma(e_1,e_2),H\rangle=0$.
Indeed, we get
$$
\begin{array}{lll}
\langle(JH)^{\perp},\nabla^{\perp}_{e_i}V\rangle
&=&\langle(JH)^{\perp},\nabla^N_{e_i}V\rangle
=-\langle\nabla^N_{e_i}(JH)^{\perp},V\rangle\\
\\&=&-\langle\nabla^N_{e_i}JH,V\rangle+\langle\nabla^N_{e_i}(JH)^{\top},V\rangle\\
\\&=&\langle JA_He_i,V\rangle+\langle\sigma(e_i,(JH)^{\top}),V\rangle\\ \\&=&0
\end{array}
$$
and
$$
\begin{array}{lll}
\langle(Je_j)^{\perp},\nabla^{\perp}_{e_i}V\rangle&=&-\langle\nabla^N_{e_i}(Je_j)^{\perp},V\rangle\\
\\&=&-\langle\nabla^N_{e_i}Je_j,V\rangle+\langle\nabla^N_{e_i}(Je_j)^{\top},V\rangle\\
\\&=&-\langle J\nabla_{e_i}e_j,V\rangle-\langle J\sigma(e_i,e_j),V\rangle+\langle\sigma(e_i,(Je_j)^{\top}),V\rangle\\
\\&=&0.
\end{array}
$$

Next, we shall prove that if a normal subbundle $S$ is orthogonal
to $L$, then so is $\nabla^{\perp}S$, i.e.
$$
\langle\sigma(e_i,e_j),\nabla^{\perp}_{e_k}V\rangle=0,\quad
\langle
J\sigma(e_i,e_j),\nabla^{\perp}_{e_k}V\rangle=0\quad\textnormal{and}\quad
\langle Je_i,\nabla^{\perp}_{e_k}V\rangle=0
$$
for any $V\in S$ and $i,j,k\in\{1,2\}$. Since we have just proved
the last property, it remains only to verify the first two of
them.

We denote
$A_{ijk}=\langle\nabla^{\perp}_{e_k}\sigma(e_i,e_j),V\rangle$ and,
since $\sigma$ is symmetric, we have $A_{ijk}=A_{jik}$. We also
obtain
$A_{ijk}=-\langle\sigma(e_i,e_j),\nabla^{\perp}_{e_k}V\rangle$,
since $V$ is orthogonal to $L$. We get
$$
\begin{array}{lll}
\langle(\nabla^{\perp}_{e_k}\sigma)(e_i,e_j),V\rangle&=&
\langle\nabla^{\perp}_{e_k}\sigma(e_i,e_j),V\rangle-\langle\sigma(\nabla_{e_k}e_i,e_j),V\rangle
-\langle\sigma(e_i,\nabla_{e_k}e_j),V\rangle\\
\\&=&\langle\nabla^{\perp}_{e_k}\sigma(e_i,e_j),V\rangle,
\end{array}
$$
and, from the Codazzi equation,
$$
\begin{array}{lll}
\langle(\nabla^{\perp}_{e_k}\sigma)(e_i,e_j),V\rangle&=&\langle(\nabla^{\perp}_{e_i}\sigma)(e_k,e_j)+(R^N(e_k,e_i)e_j)^{\perp},V\rangle\\
\\&=&\langle(\nabla^{\perp}_{e_j}\sigma)(e_k,e_i)+(R^N(e_k,e_j)e_i)^{\perp},V\rangle\\
\\&=&\langle(\nabla^{\perp}_{e_i}\sigma)(e_k,e_j),V\rangle=\langle(\nabla^{\perp}_{e_j}\sigma)(e_k,e_i),V\rangle.
\end{array}
$$
We have just proved that
$A_{ijk}=A_{kji}=A_{ikj}$.

Next, since $\nabla^{\perp}_{e_k}V$ is orthogonal to $JT\Sigma^2$
and to $JH$, it follows that the frame field $\{e_1,e_2\}$
diagonalizes $A_{\nabla^{\perp}_{e_k}V}$ and we get
$$
A_{ijk}=-\langle\sigma(e_i,e_j),\nabla^{\perp}_{e_k}V\rangle=-\langle
e_i,A_{\nabla^{\perp}_{e_k}V}e_j\rangle=0
$$
for any $i\neq j$. Hence, we have obtained that $A_{ijk}=0$ if two
indices are different from each other.

Finally, we only have to prove that $A_{iii}=0$. Indeed, we have
$$
\begin{array}{lll}
A_{iii}&=&-\langle\sigma(e_i,e_i),\nabla^{\perp}_{e_i}V\rangle
=-\langle 2H,\nabla^{\perp}_{e_i}V\rangle+\langle\sigma(e_j,e_j),\nabla^{\perp}_{e_i}V\rangle\\
\\&=&\langle 2\nabla^{\perp}_{e_i}H,V\rangle-A_{jji}=0.
\end{array}
$$
It is easy to see that if $V$ is orthogonal to $L$, then $JV$ is
normal and orthogonal to $L$. It follows that
$$
\begin{array}{lll}
\langle(J\sigma(e_i,e_j))^{\perp},\nabla^{\perp}_{e_k}
V\rangle &=&-\langle\nabla^N_{e_k}(J\sigma(e_i,e_j))^{\perp},V\rangle\\
\\&=&-\langle\nabla^N_{e_k}J\sigma(e_i,e_j),V\rangle+\langle\nabla^N_{e_k}(J\sigma(e_i,e_j))^{\top},V\rangle\\ \\
&=&\langle JA_{\sigma(e_i,e_j)}e_k,V\rangle-\langle J\nabla^{\perp}_{e_k}\sigma(e_i,e_j),V\rangle\\
\\&&+\langle\sigma(e_k,(J\sigma(e_i,e_j))^{\top}),V\rangle\\ \\&=&\langle\nabla^{\perp}_{e_k}\sigma(e_i,e_j),JV\rangle
\\ \\&=&0.
\end{array}
$$

Thus, we come to the conclusion that the subbundle $L$ is
parallel.
\end{proof}

In the case when $H$ is umbilical we can use the quadratic form
$Q$ to prove the following

\begin{lemma}\label{lemma_umbilical} Let $\Sigma^2$ be an immersed surface in a complex space form
$N^{n}(\rho)$, $\rho\neq 0$, with nonzero parallel mean curvature
vector $H$. If $H$ is an umbilical direction everywhere, then
$\Sigma^2$ is a totally real pseudo-umbilical surface of $N$.
\end{lemma}

\begin{proof} Since $H$ is umbilical it follows that
$\langle\sigma(Z,Z),H\rangle=0$, which implies that $\Sigma^2$ is
pseudo-umbilical and that $Q(Z,Z)=3\rho\langle JZ, H\rangle^2$.

Next, as the $(2,0)$-part of $Q$ is holomorphic, we have $\bar
Z(Q(Z,Z))=0$ and further
$$
0=\bar Z(\langle JZ,H\rangle^2)=-2|H|^2\langle JZ, H\rangle\langle
JZ,\bar Z\rangle,
$$
as we have seen in a previous section. Hence, $\langle JZ,\bar
Z\rangle=0$ or $\langle JZ,H\rangle=0$. Assume that the set of
zeroes of $\langle JZ,\bar Z\rangle=0$ is not the entire
$\Sigma^2$. Then, by analyticity, it is a closed set without
interior points and its complement is an open dense set in
$\Sigma^2$. In this last set we have $\langle JZ,H\rangle=0$ and
then, since $H$ is parallel and $\Sigma^2$ is pseudo-umbilical,
$$
\begin{array}{lll}
0=\bar Z(\langle JZ, H\rangle)&=&\langle J\nabla^N_{\bar Z}Z,
H\rangle+\langle JZ,\nabla^N_{\bar Z}H\rangle\\ \\&=&-\langle\bar
Z,Z\rangle\langle JH,H\rangle-\langle JZ,A_H\bar Z\rangle\\
\\&=&-|H|^2\langle JZ,\bar Z\rangle,
\end{array}
$$
which means that $\Sigma^2$ is also totally real.
\end{proof}

\begin{remark}\label{rem_CO} Some kind of a converse result was obtained by
B.-Y. Chen and K. Ogiue since they proved in \cite{CO} that if a
unit normal vector field to a $2$-sphere, immersed in a complex
space form as a totally real surface, is parallel and
isoperimetric, then it is umbilical.
\end{remark}

\begin{remark} In \cite{S} N. Sato proved that, if $M$ is a
pseudo-umbilical submanifold of a complex projective space
$\mathbb{C}P^n(\rho)$, with nonzero parallel mean curvature vector
field, then it is a totally real submanifold. Moreover, the mean
curvature vector field $H$ is orthogonal to $JTM$. Therefore, if
$M$ is a surface, it follows that the $(2,0)$-part of $Q$ vanishes
on $M$.
\end{remark}

\begin{remark} In order to show that only the two situations exposed in Lemma
\ref{lemma_parallel} and Lemma \ref{lemma_umbilical} can occur, we
shall use an argument similar to that in Remark $5$ in \cite{AdCT}.
Thus, since the map $p\in\Sigma^2\rightarrow(A_H-\mu\id)(p)$,
where $\mu$ is a constant, is analytic, it follows that if $H$ is
an umbilical direction, then this either holds on $\Sigma^2$ or
only for a closed set without interior points. In this second case
$H$ is not an umbilical direction in an open dense set, and then Lemma \ref{lemma_parallel} holds on this set. By continuity it holds on $\Sigma^2$.
\end{remark}

By using Lemma \ref{lemma_parallel} and Lemma
\ref{lemma_umbilical} we can state

\begin{proposition}\label{split} Either $H$ is everywhere an umbilical
direction, and $\Sigma^2$ is a totally real pseudo-umbilical
surface of $N$, or $H$ is nowhere an umbilical direction, and
there exists a subbundle of the normal bundle that is parallel,
contains the image of the second fundamental form and its
dimension is less or equal to $8$.
\end{proposition}

Now, from Proposition \ref{split} and a result of J. H. Eschenburg
and R. Tribuzy (Theorem 2 in \cite{ET}), it follows

\begin{theorem}\label{th_red} Let $\Sigma^2$ be an isometrically immersed surface in a complex
space form $N^n(\rho)$, $n\geq 3$, $\rho\neq 0$, with nonzero
parallel mean curvature vector. Then, one of the following holds:
\begin{enumerate}

\item $\Sigma^2$ is a totally real pseudo-umbilical surface of
$N^n(\rho)$, or

\item $\Sigma^2$ is not pseudo-umbilical and it lies in a complex
space form $N^r(\rho)$, where $r\leq 5$.
\end{enumerate}
\end{theorem}

\begin{remark} The case when $\rho=0$ is solved by Theorem 4 in \cite{Y}.
\end{remark}

\begin{remark} We have seen (Remark \ref{rem_CO}) that if $\Sigma^2$ is a totally real
$2$-sphere then it is pseudo-umbilical and therefore the second
case of the previous Theorem cannot occur for such surfaces.
\end{remark}

\section{$2$-spheres with constant K\"ahler angle in complex space forms}

This section is devoted to the study of immersed surfaces
$\Sigma^2$ in a complex space form $N^{n}(\rho)$, $n\geq 3$,
$\rho\neq 0$, with nonzero non-umbilical parallel mean curvature
vector $H$ and constant K\"ahler angle, on which the $(2,0)$-part
of $Q$ vanishes. We shall compute the Laplacian of the function
$|A_H|^2$ for such a surface and show that there are no
$2$-spheres with these properties.

Let $\{e_1,e_2\}$ be an orthonormal frame on $\Sigma^2$ such that
$H\perp Je_1$. The fact that the $(2,0)$-part of the quadratic
form $Q$ vanishes can be written as
\begin{equation}\label{eq:q=0}
\begin{cases}
8|H|^2\langle\sigma(e_1,e_1)-\sigma(e_2,e_2),H\rangle=-3\rho(\langle
Je_1,H\rangle^2-\langle Je_2,H\rangle^2)\\
8|H|^2\langle\sigma(e_1,e_2),H\rangle=3\rho\langle
Je_1,H\rangle\langle Je_2,H\rangle,
\end{cases}
\end{equation}
and, from the second equation, we see that
$\langle\sigma(e_1,e_2),H\rangle=0$. It follows that the frame
$\{e_1,e_2\}$ diagonalizes simultaneously $A_H$ and $A_V$, for all
normal vectors $V$ as in Corollary \ref{lemma_split}, since we are
in the second case of Theorem \ref{th_red}.

Next, since $\Sigma^2$ is not holomorphic or anti-holomorphic, we
have $\cos\theta\neq\pm 1$ on an open dense set and we can
consider again the normal vectors
$$
e_3=-\cot\theta
e_1-\frac{1}{\sin\theta}Je_2\quad\textnormal{and}\quad
e_4=\frac{1}{\sin\theta}Je_1-\cot\theta e_2
$$
and obtain an orthonormal frame $\{e_1,e_2,e_3,e_4\}$ in
$\Span\{e_1,e_2,Je_1,Je_2\}$, where $\theta$ is the K\" ahler
angle on $\Sigma^2$.

It is easy to see that if $H\perp JT\Sigma^2$ it results that the
surface is pseudo-umbilical, which is a contradiction.

On the other hand, if we assume that $H\in\Span\{e_3,e_4\}$ it
follows $H=\pm|H|e_3$, since $Je_1\perp H$, and then $e_3$ is
parallel. Also, since all normal vectors but $e_4$ verify
conditions in Corollary \ref{lemma_split} we have
$\sigma(e_1,e_2)\parallel e_4$. By using these facts and the
expression of $e_3$ we obtain that
$\sigma(e_i,e_j)\in\Span\{e_3,e_4\}$ for $i,j\in\{1,2\}$, and then
$\dim L=2$, where $L$ is the subbundle in Lemma
\ref{lemma_parallel}. Therefore, again by the meaning of Theorem 2
in \cite{ET}, we get that $\Sigma^2$ lies in a complex space form
$N^2(\rho)$, which case was studied earlier in this paper.

Consequently, in the following, we shall assume that
$H\notin\Span\{e_3,e_4\}$, and, as we also know that $H$ is not
orthogonal to $JT\Sigma^2$, it results that the mean curvature
vector can be written as
$$
H=|H|(\cos\beta e_3+\sin\beta e_5)
$$
where $\beta$ is a real-valued function defined locally on
$\Sigma^2$ and $e_5$ is a unit normal vector field such that
$e_5\perp JT\Sigma^2$. We consider the orthonormal frame field
$$
\{e_1,e_2,e_3,e_4,e_5,e_6=Je_5,\ldots,e_{2n-1},e_{2n}=Je_{2n-1}\}
$$
on $N$ and its dual frame $\{\theta_i\}_{i=1}^{2n}$. These are
well defined at the points of $\Sigma^2$ where $\sin(2\beta)\neq
0$, which, due to our assumptions, form an open dense set in
$\Sigma^2$. The structure equations of the surface are
$$
d\phi=-i\theta_{12}\wedge\phi\quad\textnormal{and}\quad
d\theta_{12}=-\frac{i}{2}K\phi\wedge\bar\phi,
$$
where $\phi=\theta_1+i\theta_2$, the real $1$-form $\theta_{12}$
is the connection form of the Riemannian metric on $\Sigma^2$ and
$K$ is the Gaussian curvature.

A result of T. Ogata in \cite{O2}, together with $H\perp e_i$ for
any $i\geq 4$, $i\neq 5$, imply that, w.r.t. the above orthonormal
frame, the components of the second fundamental form are
\begin{footnotesize}
$$
\sigma^3=\left(\begin{array}{cc}|H|\cos\beta-\Re(\bar a+c)&-\Im(\bar a+c)\\
\\-\Im(\bar a+c)&|H|\cos\beta+\Re(\bar
a+c)\end{array}\right),\quad
\sigma^4=\left(\begin{array}{cc}\Im(\bar a-c)&-\Re(\bar a-c)\\
\\-\Re(\bar a-c)&-\Im(\bar a-c)\end{array}\right)
$$
$$
\sigma^5=\left(\begin{array}{cc}|H|\sin\beta-\Re(\bar a_3+c_3)&-\Im(\bar a_3+c_3)\\
\\-\Im(\bar a_3+c_3)&|H|\sin\beta+\Re(\bar
a_3+c_3)\end{array}\right)
$$
$$
\sigma^6=\left(\begin{array}{cc}\Im(\bar a_3-c_3)&-\Re(\bar a_3-c_3)\\
\\-\Re(\bar a_3-c_3)&-\Im(\bar a_3-c_3)\end{array}\right)
$$
$$
\sigma^{2\alpha-1}=\left(\begin{array}{cc}-\Re(\bar a_{\alpha}+c_{\alpha})&-\Im(\bar a_{\alpha}+c_{\alpha})\\
\\-\Im(\bar a_{\alpha}+c_{\alpha})&\Re(\bar
a_{\alpha}+c_{\alpha})\end{array}\right),\quad
\sigma^{2\alpha}=\left(\begin{array}{cc}\Im(\bar a_{\alpha}-c_{\alpha})&-\Re(\bar a_{\alpha}-c_{\alpha})\\
\\-\Re(\bar a_{\alpha}-c_{\alpha})&-\Im(\bar a_{\alpha}-c_{\alpha})\end{array}\right)
$$
\end{footnotesize}where $a,c,a_{\alpha},c_{\alpha}$, with $\alpha\in\{3,\ldots,n\}$, are
complex-valued functions defined locally on the surface
$\Sigma^2$. We note that, since $\sigma(e_1,e_2)\perp H$ and
$\sigma(e_1,e_2)\perp e_5$, it follows $\sigma(e_1,e_2)\perp e_3$.
Moreover, since $\sigma(e_1,e_2)\perp e_i$ for any
$i\in\{1,\ldots,2n\}\setminus\{4,6\}$, we have $\bar
a+c\in\mathbb{R}$, $\bar a_3+c_3\in\mathbb{R}$ and
$a_{\alpha}=c_{\alpha}$ for any $\alpha\geq 4$.

In the same paper \cite{O2}, amongst others, the author computed
the differential of the K\" ahler angle function $\theta$ for a
minimal surface. In the same way, this time for our surface, we
get
$$
d\theta=\Big(a-\frac{|H|}{2}\cos\beta\Big)\phi+\Big(\bar
a-\frac{|H|}{2}\cos\beta\Big)\bar\phi.
$$

The next step is to determine the connection form $\theta_{12}$
and the differential of the function $\beta$, by using the
property of $H$ being parallel. We have
\begin{equation}\label{eq:H}
\nabla^{\perp}_{e_i}H=(-\sin\beta e_3+\cos\beta
e_5)d\beta(e_i)+\cos\beta\nabla^{\perp}_{e_i}e_3+\sin\beta\nabla^{\perp}_{e_i}e_5=0
\end{equation}
for $i\in\{1,2\}$, and then
$$
\cos\beta\langle\nabla^N_{e_i}e_3,e_4\rangle+\sin\beta\langle\nabla^N_{e_i}e_3,e_4\rangle=0,\quad
i\in\{1,2\}
$$
from where, by using the expressions of $e_3$ in the first term,
of $e_4$ in the second one and of the second fundamental form of
$\Sigma^2$, we get
$$
\theta_{12}(e_1)=\cot\theta\Im(\bar
a-c)-\frac{\tan\beta}{\sin\theta}\Im(\bar a_3-c_3)
$$
$$
\theta_{12}(e_2)=-|H|\frac{\cot\theta}{\cos\beta}-2\cot\theta\Re
a+\tan\beta\Big(\tan\Big(\frac{\theta}{2}\Big)\Re
a_3-\cot\Big(\frac{\theta}{2}\Big)\Re c_3\Big)
$$
and finally $\theta_{12}=f_1\phi+\bar f_1\bar\phi$, where
\begin{equation}\label{eq:f1}
f_1=\frac{i}{2}\Big(|H|\frac{\cot\theta}{\cos\beta}+2\cot\theta
a-\frac{\tan\beta}{\sin\theta}(a_3-\bar
c_3)+\cot\theta\tan\beta(a_3+\bar c_3)\Big).
\end{equation}

Now, from equation \eqref{eq:H}, we also obtain
$$
d\beta(e_i)+\langle\nabla^N_{e_i}e_3,e_5\rangle=0,\quad
i\in\{1,2\}
$$
and then, replacing $e_3$ with its expression and also using the
expression of the second fundamental form, we get
$$
d\beta(e_1)=|H|\cot\theta\sin\beta+\tan\Big(\frac{\theta}{2}\Big)\Re
a_3-\cot\Big(\frac{\theta}{2}\Big)\Re c_3,\quad
d\beta(e_2)=\frac{1}{\sin\theta}\Im(\bar a_3-c_3).
$$
Hence the differential of $\beta$ is given by $d\beta=f_2\phi+\bar
f_2\bar \phi$, where
\begin{equation}\label{eq:f2}
f_2=\frac{1}{2}\Big(|H|\cot\theta\sin\beta+\frac{1}{\sin\theta}(a_3-\bar
c_3)-\cot\theta(a_3+\bar c_3)\Big).
\end{equation}
We note that if the K\" ahler angle $\theta$ is constant, then
$a=\bar a=\frac{|H|}{2}\cos\beta$, and, from \eqref{eq:f1}, it
results
\begin{equation}\label{eq:f1const}
f_1=\frac{i}{2}\Big\{|H|\cot\theta\Big(\cos\beta+\frac{1}{\cos\beta}\Big)-\frac{\tan\beta}{\sin\theta}(a_3-\bar
c_3)+\cot\theta\tan\beta(a_3+\bar c_3)\Big\}.
\end{equation}

Let us now return to the first equation of \eqref{eq:q=0}, which
can be rewritten as
$$
\mu_1-\mu_2=\frac{3}{8}\rho\sin^2\theta\cos^2\beta,
$$
where $A_He_i=\mu_ie_i$. Since $\mu_1+\mu_2=2|H|^2$ we have
$\mu_1=|H|^2+\frac{3}{16}\rho\sin^2\theta\cos^2\beta$ and
$\mu_2=|H|^2-\frac{3}{16}\rho\sin^2\theta\cos^2\beta$. Thus
\begin{equation}\label{eq:Abound}
|A_H|^2=\mu_1^2+\mu_2^2=2|H|^4+\frac{9}{128}\rho^2\sin^4\theta\cos^4\beta.
\end{equation}

In the following, we shall assume that the K\" ahler angle of the
surface $\Sigma^2$ is constant and then the Laplacian of $|A_H|^2$
is given by
$$
\Delta|A_H|^2=\frac{9}{128}\rho^2\sin^4\theta\Delta(\cos^4\beta).
$$

In order to compute the Laplacian of $\cos^4\beta$ we need the
following formula, obtained by using \eqref{eq:f2} and
\eqref{eq:f1const},
$$
\begin{array}{ll}
d(\cos^4\beta)&=-4\sin\beta\cos^3\beta
d\beta=-4\sin\beta\cos^3\beta(f_2\phi+\bar f_2\bar\phi)\\
\\&=-4\cos^4\beta\Big\{\Big(if_1+|H|\frac{\cot\theta}{\cos\beta}\Big)\phi+
\Big(-i\bar
f_1+|H|\frac{\cot\theta}{\cos\beta}\Big)\bar\phi\Big\}.
\end{array}
$$
We also have
$dd^c(\cos^4\beta)=\frac{i}{2}(\Delta(\cos^4\beta))\phi\wedge\bar\phi$
and
$$
d^c(\cos^4\beta)=-4i\cos^4\beta\Big\{\Big(-i\bar
f_1+|H|\frac{\cot\theta}{\cos\beta}\Big)\bar\phi-\Big(if_1+|H|\frac{\cot\theta}{\cos\beta}\Big)\phi\Big\}.
$$
After a straightforward computation, we get
$$
\Delta(\cos^4\beta)=4\cos^4\beta\Big(K+4|f_1|^2+12\Big|if_1+|H|\frac{\cot\theta}{\cos\beta}\Big|^2\Big)
$$
and then
$$
\Delta|A_H|^2=\frac{9}{32}\rho^2\sin^4\theta\cos^4\beta\Big(K+4|f_1|^2+12\Big|if_1+|H|\frac{\cot\theta}{\cos\beta}\Big|^2\Big).
$$

Assume now that the surface $\Sigma^2$ is complete and it has
nonnegative Gaussian curvature. It follows, from a result of A.
Huber in \cite{AH}, that $\Sigma^2$ is parabolic. Then, from the
above formula, we get that $|A_H|^2$ is a subharmonic function,
and, since $|A_H|^2$ is bounded (due to \eqref{eq:Abound}), it
results $K=0$, which, together with the Gauss-Bonnet Theorem, lead
to the following non-existence result.

\begin{theorem} There are no $2$-spheres with nonzero non-umbilical
parallel curvature vector and constant K\"ahler angle in a
non-flat complex space form.
\end{theorem}

\end{document}